\title{Conflict-free coloring of graphs}
\author{{Roman Glebov \thanks{
Department of Mathematics, ETH, 8092 Zurich, Switzerland.
Previous affiliation: Mathematics Institute and DIMAP, University of Warwick, Coventry CV4 7AL, UK.
This research was done when the author was affiliated with the
Institute of Mathematics,
Free University Berlin, 14195 Berlin, Germany. Email:
roman.l.glebov@gmail.com.
The author was supported by DFG within the research training group ``Methods for Discrete Structures''.}}
\quad {Tibor Szab\'o
\thanks{ Institute of Mathematics,
Free University Berlin, 14195 Berlin, Germany, email:
szabo@math.fu-berlin.de.}}
\quad {G\'abor Tardos \thanks{Alfr\'ed R\'enyi Institute of Mathematics,
            Hungarian Academy of Sciences,
            Budapest, Hungary and Zhejiang Normal
            University, Jinhua, China, email: tardos@renyi.hu.
            Partially supported by an NSERC grant, the MTA Cryptography
            ``Lend\"ulet'' project, and the Hungarian OTKA grant NN-102029.}}
}

\documentclass[headsepline,a4paper,openany]{article}
\usepackage{amsthm,amsmath,amssymb,amsfonts}
\usepackage{epsfig}
\usepackage{graphicx,color}
\usepackage{pstricks}
\usepackage{verbatim}
\usepackage[T1]{fontenc}
\usepackage[ansinew]{inputenc}
\usepackage{stmaryrd}
\usepackage[colorlinks]{hyperref}

\newtheorem{lem}{Lemma}

\newtheorem{teo}{Theorem}

\newtheorem{obs}{Observation}

\newcommand{\proofstart}{{\bf Proof\hspace{2em}}}

\newcommand{\proofend}{\hspace*{\fill}\mbox{$\Box$}\vspace*{7mm}}
\let\epsilon=\varepsilon
\def\probnew{\mathbb{P}}

\begin{document}

\maketitle

\begin{abstract}
We study the conflict-free chromatic number $\chi_{CF}$
of graphs from extremal and probabilistic point of view.
We resolve a question of Pach and Tardos about the
maximum conflict-free chromatic number an $n$-vertex graph can
have. Our construction is randomized. In relation to this
we study the evolution of the conflict-free chromatic number
of the Erd\H os-R\'enyi random graph $G(n,p)$ and give
the asymptotics for $p=\omega(1/n)$.
We also show that for $p\geq 1/2$ the conflict-free chromatic
number differs from the domination number by at most $3$.
\end{abstract}

{\bf MSC classes:} 05C35, 05C15, 05C80, 05D40, 05C69.

\section{Introduction and definitions}

Let $G=(V,E)$ be a simple graph. For every $x\in V$ we denote by $N(x)=\{y\in V:~xy\in E\}$ its {\em neighborhood}
and by $N[x]= N(x)\cup \{ x\}$ its {\em closed neighborhood}.
A (not necessarily proper) vertex coloring $\chi$ of $G$ is called {\em conflict-free},
if for each vertex $x\in V$, there exists a vertex $y$ in $N[x]$ whose color is different from
the color of each other vertex in $N[x]$.
We then say that $y$ has {\em unique
color} in $N[x]$.
The {\em conflict-free chromatic number} $\chi_{CF}(G)$ is
the smallest $r$, such that there exists a conflict-free $r$-coloring
of $G$.
Conflict-free coloring can be interpreted as a relaxation of the usual
{\em proper coloring} concept where each vertex $x$ is required to have
a unique color in its
own closed neighborhood $N[x]$. Hence $\chi_{CF}(G) \leq \chi(G)$ for every
graph $G$.

The study of conflict-free colorings was originated in the work of Even, Lotker, Ron, and Smorodinsky~\cite{even} and
Smorodinsky~\cite{shakharsdiss} who were motivated by the problem of frequency assignment in cellular
networks.
(See the recent survey by Smorodinsky~\cite{shakhar}.)
In most of these classical instances the graphs studied arise
from a geometric setting.
Recently Pach and Tardos~\cite{PT09} initiated the study of the problem for abstract graphs and hypergraphs.
Here we continue the consideration of conflict-free colorings of abstract graphs.

Note that, unlike the proper coloring number, the conflict-free chromatic number is not monotone.
In particular, in the two extremes $\chi_{CF}(K_n)=2$ for the complete graph and
$\chi_{CF}(\bar{K}_n)=1$ for the empty graph, while the conflict-free
chromatic number of general graphs can be arbitrarily high. We investigate this parameter from extremal and
probabilistic points of view.

Pach and Tardos~\cite{PT09} raised the problem of determining the order of
magnitude of $\chi_{CF}(n):=\max \{\chi_{CF}(G):~|V(G)|=n\}$, the largest
conflict-free chromatic number an $n$-vertex graph can have.
From above they showed $\chi_{CF}(n) = O\left(\ln^2 n\right)$ but from below they could
only prove that the conflict-free coloring number of
the random graph $G\left(n,\frac{1}{2}\right)$
is asymptotically almost surely $\Omega(\ln n)$, hence
$\chi_{CF}(n) = \Omega(\ln n)$. Here asymptotically almost surely means
probabilities tending to 1 as $n$ goes to infinity and it will be abbreviated
below as a.a.s.

At first one could try to improve the lower bound $\chi_{CF}(n)$ by
considering the random graph $G(n,p)$
with some $p=p(n)\ne1/2$.
In our first theorem we give tight
estimates (holding a.a.s.) for the conflict-free chromatic number of these
random graphs. Our bounds show that some probabilities $p(n)\rightarrow 0$
yield the highest conflict-free coloring numbers for the $G(n,p(n))$, but
these are only a constant factor larger than those of $G(n,1/2)$.

To state our theorem we introduce
\[\mu=\mu(p)=\max \{ip(1-p)^{i-1}:~i\in \mathbb{N}^+\}\]
for $0 < p  < 1$. Notice that the maximum is taken at $i=\lfloor1/p\rfloor$ so
we have
\[\mu(p)=\left\lfloor \frac{1}{p}\right\rfloor p(1-p)^{\lfloor 1/p\rfloor
  -1},\]
and (as simple calculation shows) this is a strictly increasing function tending to
$e^{-1}$ as $p$ goes to $0$.

\begin{teo}
\label{random-graph}
For every $\epsilon>0$ and every function $0<p=p(n)<1-\epsilon$ such that $np(n)\rightarrow\infty$, the
following holds a.a.s.\
\[ (1-\epsilon)\frac{\ln (np)}{-\ln (1-\mu(p))}\leq
\chi_{CF}(G(n,p))\leq (1+\epsilon)\frac{\ln (np)}{-\ln (1-\mu(p))}.\]
\end{teo}
Note that the theorem implies $\chi_{CF}(G(n,p))=O(\log n)$ a.a.s.
for all $p$ considered. It is not hard to show that the $O(\log n)$ upper bound is also valid a.a.s. in the full range of $p\in [0,1]$.
For the range $p = O(1/p)$ this follows from
$\chi_{CF}(G) \leq \chi(G) \leq \Delta(G)+1$, but in this range we are not able to determine the asymptotics.

For $1/2\le p< 1$ we can prove an even
tighter result: the conflict-free
coloring number differs by at most 3 from
the domination number.
A set $S$ of vertices of a graph $G$ constitutes a
{\em dominating set} if each $v \in V$ is
either in $S$ or is adjacent to a vertex in $S$.
The {\em domination number} $D(G)$ is the smallest size of
a dominating set in $G$.
\begin{teo}
\label{random-graph-domination} For every graph $G$,
\[ \chi_{CF}(G)\leq  D(G)+ 1.\]
Furthermore, for $\frac{1}{2}\le p(n)$ a.a.s.\
\[D(G(n,p(n)))-3 \leq \chi_{CF}(G(n,p(n))).\]
\end{teo}
The domination number of the random graph with constant $p$ was pinned down to
be one of two integers a.a.s.\ by Wieland and Godbole~\cite{WG01}.
Furthermore, it was observed by Glebov, Liebenau, and Szab\'o~\cite{ichanitatibor} that the same result
holds also for a variable $p(n)$. The following is a corollary of these results for the range of our interest:

\begin{teo}[Corollary of \cite{WG01,ichanitatibor}]\label{domina}
For $1/2\le p<1$ the domination number $D(G(n,p(n))$ is either
$\left\lfloor\frac{\ln n-2\ln\ln n+\ln\ln
    \frac{1}{1-p}}{-\ln(1-p)}\right\rfloor + 1$ or one more a.a.s.
\end{teo}

Hence the behavior of $\chi_{CF}$ is also very well understood in this range.
In fact, we prove Theorem~\ref{random-graph-domination} by calculating the
a.a.s.\ lower bound on $\chi_{CF}(G(n,p))$
and comparing it with the a.a.s.\ domination number.
Notice that using Theorem~\ref{domina}
Theorem~\ref{random-graph-domination} implies Theorem~\ref{random-graph} for
the range $p\geq 1/2$, where we have $\mu(p)=p$. We mention that for
the range $p<\frac{1}{2}$ the results
of~\cite{WG01,ichanitatibor} on the domination number and
Theorem~\ref{random-graph} imply that the conflict-free
chromatic number and the domination
number {\em differ} in the asymptotics.

In our final result we resolve the open problem of Pach and Tardos~\cite{PT09} regarding
$\chi_{CF} (n)$ by constructing $n$-vertex graphs
$G$ with $\chi_{CF}(G) = \Omega(\ln^2 n)$.

\begin{teo}
\label{main}
\[\chi_{CF}(n) = \Theta(\ln^2 n).\]
\end{teo}

The structure of the paper is the following:
in Section~\ref{sec:random-graph} we prove Theorems~\ref{random-graph}
and \ref{random-graph-domination},
while Theorem~\ref{main} is proven in Section~\ref{sec:main}.
For simplicity we routinely omit floor and ceiling signs
as long as they do not influence the validity of our asymptotic statements.

{\bf Notation.}
Let $G$ be a graph with vertex set $V=V(G)$ and let $A\subseteq V$. We say that $N_G^{(1)}(A)=\{v \in V\setminus A: ~
|N(v)\cap A|=1\}$ is the {\em one-neighborhood of $A$} and
$\overline{N_G}(A)=V\setminus\bigcup_{x\in A}N[x]$ is the {\em
non-neighborhood of $A$}. The subscript $G$ is omitted if it is clear from the
context.

We use $V\choose m$ to denote the set of all $m$-element subsets of $V$.

\section{Evolution of the conflict-free chromatic number
in random graphs}
\label{sec:random-graph}

\subsection{Upper bounds}

A simple upper bound is obtained from the fact that any proper coloring is a conflict-free coloring, so
\[\chi_{CF}(G) \leq \chi(G).\]
However, this bound is a.a.s.\ not tight for the random graph $G(n,p)$ in the range of $p$ we are interested in,
i.e., for $p=\omega(1/n)$.

Another inequality involves domination.
If a set of vertices $S$ is a dominating set of
$G$ then one can construct a conflict-free coloring of $G$
with $|S|+1$ colors by
giving $|S|$ distinct colors to the vertices in $S$ and
one further color to vertices in $V(G)\setminus S$.
Hence for every graph $G$
\[\chi_{CF}(G) \leq D(G)+1.\]
This proves the upper bound in Theorem~\ref{random-graph-domination}.

The rest of this section deals with the upper bound in
Theorem~\ref{random-graph}.

Regarding conflict-free colorings
the crucial property of a vertex $x$ is whether it has exactly one neighbor
in some color class $S$ and hence the color of $S$ is unique in $N[x]$.
For a fixed set $S$ and a fixed vertex $x\in V\setminus S$
the probability of this happening is $|S|p(1-p)^{|S|-1}$.
This motivates our definition of $\mu(p)$ in Section~1 as the maximum of this
probability for any color class size. We let $m=\lfloor1/p\rfloor$ stand for
the ``most desirable'' color class size maximizing the above probability and
giving $\mu =mp(1-p)^{m-1}$.

Since the upper bound in Theorem~\ref{random-graph-domination} implies the
upper bound in Theorem~\ref{random-graph}
for $p\geq\frac{1}{2}$, we assume $p<\frac{1}{2}$ from now on.
To start we prove two technical lemmas for random graphs.

First we give an explicit bound on the probability that the domination number
of a random graph is extremely low. We need the explicit bound in because we
will use the union bound for more than a constant number of similar events and
thus the a.a.s.\ bound of Theorem~\ref{domina} is not enough.

\begin{lem}
\label{l1}
For any $\ell\in\mathbb N^+$ and $p$ with $ \frac{100}{\ell}<p<\frac{1}{2}$
we have that
\[ \probnew\left[D(G(\ell,p))<m\right]< 0.9^\ell,\]
where $m=\lfloor1/p\rfloor$ as before.
\end{lem}
\proofstart
Throughout the proof we will use that
$m-1 < \frac{1}{p} < \frac{\ell}{100}$.
Let $S\subset V$ be a set of size $m-1$. The probability that a vertex
$x \in V \setminus S$ has no neighbor in $S$ is
\[
\probnew[N(x)\cap S = \emptyset] = (1-p)^{|S|} \geq (1-p)^{1/p}> 1/4.
\]
The events that $N(x)\cap S = \emptyset$ are independent  for $x\in V\setminus S$, hence $S$ is dominating with
probability $<(3/4)^{\ell-m+1}$.
The probability in the lemma is
\begin{align*}
\probnew \left[\exists S \in {V \choose m-1}: \overline{N}(S) = \emptyset \right]
&< {\ell\choose m-1} (3/4)^{\ell-m+1}\\
&< 0.9^\ell.
\end{align*}
\proofend

We know that the expected size of the one-neighborhood of a set of
vertices of size $m$ is $(|V|-m)\mu$. The following is a routine observation that
the actual size deviates largely from this expectation with a very low
probability.
\begin{lem}
\label{l2}
For every $\delta>0$ there exists a $K=K(\delta)$ such that for any $p=p(\ell)>
\frac{K}{\ell}$ in $G(\ell,p)$
we have that
\[\probnew\left[\exists S \in {V \choose m}: |N^{(1)}(S)| < (1-\delta)\mu (\ell-m)\right ]
< e^{-\frac{\delta^2}{4}\mu\ell },\]
where $m=\lfloor1/p\rfloor$ and $\mu=\mu(p)=mp(1-p)^{m-1}$.
\end{lem}
\proofstart
For an arbitrary set $S\subset V$ of size $m$ and vertex $x \in V \setminus S$,
the probability that $x$ has exactly one neighbor in $S$ is $\mu$.
The random variable $\left|N^{(1)}(S)\right|$ is the sum of $\ell-m$ mutually
independent characteristic variables and its expectation is $\mu(\ell-m)$.
Hence by the Chernoff bound and the union bound we have
\begin{align*}
\probnew\left[\exists S \in {V \choose m }: |N^{(1)}(S)| < (1-\delta)\mu(\ell-m) \right]
&< {\ell\choose m } e^{-\frac{\delta^2}{2}\mu (\ell-m)}\\
&\leq \left( \left(eK\right)^{\frac{1}{K}}
e^{-\frac{\delta^2}{3}\mu }\right)^\ell,
\end{align*}
and the bound follows if $K$ is sufficiently large.
\proofend

Let us choose $\delta=\delta(\epsilon)>0$ such that it satisfies
\begin{align}
\frac{1+\epsilon}{-\ln(1-\mu(p))}>\frac1{-\ln(1-(1-\delta)\mu(p))}+\delta\label{definequal}
\end{align}
and assume $K=K(\delta)$ from the Lemma~\ref{l2} satisfies $K>100$ so we can
also use Lemma~\ref{l1}. Assuming $p=p(n)$ satisfies $np\rightarrow\infty$
(or even the weaker condition $p>\frac{K^*}{n}$ for $K^*=e^{K/\delta}$) we
give a deterministic
algorithm which a.a.s.\ constructs a conflict-free coloring of $G(n,p)$
using $(1+\epsilon)\frac{\ln (np)}{-\ln(1-\mu)}$ colors. In this algorithm
$d(G_i)$ denotes the {\em degeneracy} of the graph $G_i$, i.e., the largest
minimum degree a non-empty subgraph of $G_i$ has.

\begin{figure}[htb]
\begin{tabbing}
....\=....\=....\=....\=.................\kill\\
{\bf \underline{Algorithm CFC}}($G$, $p$, $\delta$)\\
\textbf{Input:} graph $G$, $V(G)=[n]$, $p\in [0,1]$, $\delta>0$.\\
Set $G_1:=G$, $n_1:=n$, $i:=1$, $m=\left \lfloor \frac{1}{p}\right\rfloor$, $\mu=mp(1-p)^{m-1}$, $K=K(\delta)>100$.\\
\textbf{while} $n_i > \ln\ln n$ and
$p > \frac{K}{n_i }$ \textbf{do} \\
\> select an independent set $S_i$ by starting with $S_i=\emptyset$ and iteratively adding\\
\> the smallest vertex in $\overline{N_{G_i}}(S_i)$
until either $\overline{N_{G_i}}(S_i)=\emptyset$ or $|S_i|=m$.\\
\>Color vertices in $S_i$ with color $i$,
color vertices in $N^{(1)}_{G_i}(S_i)$ with color $0$,\\
\>define $G_{i+1}:=G_i-\left(S_i \cup N^{(1)}_{G_i}\left(S_i\right)\right)$,
$n_{i+1}:=|V(G_{i+1})|$, $i:=i+1$\\
Color $G_i$ properly using $d(G_i)+1$ new colors.
\end{tabbing}
\end{figure}

Notice that all executions of the main {\bf while}-loop of the algorithm use a
separate color and only color $0$ is used in many executions. Note also that
this color $0$ is a ``filler color'' as
it is never used as the unique color in the closed neighborhood of some vertex
to ensure the conflict-free property of the coloring is obtained.

Let $I$ be the last value of the index $i$ in the algorithm. Clearly the
algorithm colors all vertices with $I+d(G_I)+1$ colors. To see that
this coloring is conflict-free let
$w\in V(G)$ be an arbitrary vertex and let $i$ be the largest index with $w\in
V(G_i)$. If $i<I$, then there is a unique vertex in $N[w]$ of color $i$ (which
may or may not be $w$ itself). If $i=I$, then $w$ has unique color in
$N[w]$.

To finish the proof it is enough to bound the values of $I$ and
$d(G_I)$ a.a.s. We start with $I$. Note that for any $1\le i\le I$ the sets
$S_1,\ldots, S_{i-1}$ selected by the algorithm, and hence the vertex set
$V(G_i)$ as well, depend only on the edges incident to $S_1 \cup \cdots \cup
S_{i-1}$. Thus, given any way the main {\bf while}-loop is executed for the
first $i-1$ times, the graph $G_i$ is still a random graph $G(n_i,p)$.
Now we
estimate the probability that $\left|N_{G_i}^{(1)}(S_i)\right|<(1-\delta)\mu(n_i-m)$. This
can happen either with $|S_i|=m$ or with $|S_i|<m$. The probability of the
former is bounded by Lemma~\ref{l2}, while the latter implies that $S_i$ is
dominating in $G_i$, the probability of which is bounded by
Lemma~\ref{l1}. Using the explicit bounds in the lemmas and the fact that the
sizes of the graphs considered are decreasing and lower bounded by a super-constant function of $n$
we conclude that a.a.s.\ in no iteration do we have
either of these anomalies:
\[\sum_{i=1}^{I-1}e^{-\delta^2\mu n_i/4}+0.9^{n_i}\leq
\sum_{\ell=\ln\ln n}^{n}e^{-\delta^2\mu \ell/4}+0.9^{\ell}=o(1).\]
Thus a.a.s.\ we must have
$n_{i+1}\le(1-(1-\delta)\mu)n_i$ for each $i<I$. Using $n_1=n$ and
$n_{I-1}>K/p$ we have a.a.s.\
\[I<\frac{\ln(np/K)}{-\ln(1-(1-\delta)\mu)}+2.\]

It remains to show that a.a.s.\ $d(G_I)<\delta\ln(np)$ and using the defining
inequality~\eqref{definequal} for $\delta$ the upper bound in
Theorem~\ref{random-graph} follows. We use again the
observation that independent of the executions of the {\bf while}-loop
$G_I$ is a random graph $G(n_I,p)$ with $n_I$ being small enough to
trigger one of the halting conditions.

If we have $p\le K/n_I$, then the expected
degree of any vertex in $G(n_I,p)$ is less than $K$.
Hence either a.a.s.\ $d(G_I)\leq K$ by
the results of~Pittel, Spencer, and Wormald~\cite{PSW} and \L uczak~\cite{Luczak} and we are done as
$K<\delta\ln(np)$,
or $n_I$ is bounded by a constant, in which case we can color $G_I$ with $n_I\le \delta\ln(np)$ colors.
If, however, $p>K/n_I$ we must have $n_I\le\ln\ln n$ to
halt the {\bf while}-loop, so we have $\ln(np)>\ln(Kn/n_I)=\Omega(\ln
n)$. Thus we have $d(G_I)<n_I<\delta\ln(np)$ if $n$ is large enough.

Note that the $\ln\ln n$ bound in the halting condition of the algorithm can
be replaced by any function that tends to infinity and is $o(\ln n)$.

Furthermore, observe that if $d(G_I)<\delta\ln (np)$ (which happens a.a.s.),
one can efficiently color $G_I$ with $\delta\ln (np)$ colors properly.
A linear time algorithm for coloring $G_I$ with at most $d(G_I)+1$ colors first iteratively removes a lowest degree vertex from the graph, then colors them greedily in the reverse order.
Hence, the running time of the algorithm CFC is linear in the size of the input graph.

\subsection{Lower bounds}
Tardos and Pach~\cite{PT09} used the concept of
{\em universality} to show the lower
bound $\chi_{CF}\left(G\left(n,\frac{1}{2}\right)\right)= \Omega(\ln n)$. A graph $G$ is called
$k$-{\em universal} if for all sets $B\subseteq A\subseteq V(G)$ with $|A|\le k$
there exists a vertex $x\in V(G)\setminus A$ with $N(x)\cap A=B$.
We introduce a similar concept, which is more closely related to the
idea of conflict-free coloring.
We call a graph $G$ {\em $(k,f)$-spoiling}, if for any $k$
disjoint subsets $A_1, \ldots, A_k\subseteq V(G)$ with
$|A_i|\leq f$ for every $i\in [k]$,
there exists a vertex $x\in V(G)\setminus \bigcup_i A_i$ such that
for each $A_i$ we have $|N(x)\cap A_i|\neq 1$, and
for each $A_i$ with $|A_i|=f$, $|N(x)\cap A_i|\geq 2$.
The vertex $x$ is called a {\em $f$-spoiler} for $(A_1, \ldots , A_k)$
and we say that $(A_1, \ldots , A_k)$ is {\em spoiled} by $x$.
We call a graph {\em $k$-spoiling}, if it is $(k,f)$-spoiling for
{\em some} $f$.

The following observation just serves to give an intuition for the concept.
\begin{obs}
\label{pres-universality}
A $2k$-universal graph $G$ is $(k,2)$-spoiling and consequently $k$-spoiling.
\end{obs}

The next lemma is the essence of all lower bounds in
Theorem~\ref{random-graph}.
\begin{lem}
\label{pres-cf}
If $G$ is $k$-spoiling, then $\chi_{CF}(G)>k$.
\end{lem}
\proofstart
Let $G$ be $(k,f)$-spoiling for some $f$ and consider an arbitrary
$k$-coloring $\chi$ of $V(G)$. We need to show that it is not conflict-free.
We define subsets $A_1,\ldots , A_k \subseteq V(G)$.
For each color $i$ which is used less than $f$ times by $\chi$,
we define $A_i$ to be the whole color class $\chi^{-1}(\{i\})$.
For each color $i$ which is used on at least $f$ vertices by $\chi$,
we set an arbitrary $f$-subset of vertices with color $i$ to be $A_i$.
Since $G$ is $(k,f)$-spoiling we find a vertex $x$ which is an $f$-spoiler
for these sets. Clearly, $N[x]$ has no unique color, showing that $\chi$ is
not conflict-free. \proofend

We first prove the tight lower bound of Theorem~\ref{random-graph-domination}
via studying the spoilers of $G(n,p)$ for $p \geq \frac{1}{2}$.
Comparing the bound of Theorem~\ref{domina} of Wieland and Godbole~\cite{WG01}
and Glebov, Liebenau, and Szab\'o~\cite{ichanitatibor}
with the bound in the following lemma finishes the proof.
\begin{lem}
\label{concentration}
The graph $G(n,p)$ with $ 1/2\leq p<1$ is a.a.s.\ $k$-spoiling for
$k = \left\lfloor\frac{\ln n-2\ln\ln n+\ln\ln \frac{1}{1-p}-\ln 3}{-\ln(1-p)}\right\rfloor$.
\end{lem}
\proofstart
We show that $G(n,p)$ is a.a.s.\ $(k,3)$-spoiling.
Take any set $A\subseteq V$ with $|A|\le3$ and $x\in V\setminus A$. For
$A=\emptyset$ we cannot have $|N(x)\cap A|=1$. For $|A|=1$ we have
\[\probnew[|N(x)\cap A|\ne1]=1-p,\]
for $|A|=2$ we have
\[
\probnew[|N(x)\cap A | \neq 1]= 1 - 2p(1-p) \geq 1-p,
\]
and finally for $|A|= 3$ we have
\[
\probnew[|N(x)\cap A | \geq 2] = 1 - 3p(1-p)^2 - (1-p)^3 \geq 1-p.
\]

Then for any family ${\mathcal A}= \{ A_1, \ldots , A_k\}$ of $k$ sets
of size at most $f=3$, the probability that
a fixed vertex $x\in V\setminus \bigcup_{i=1}^k A_i$ is a spoiler is
\[
\probnew[x \mbox{ is a spoiler for } \mathcal{A}]\geq (1-p)^k.
\]
Thus
\[
\probnew\left[\mathcal{A} \mbox{ is not spoiled by any } x\in V\setminus
\bigcup_{i=1}^k A_i\right]\leq \left(1-(1-p)^k\right)^{n-3k}.
\]
There are at most $(n+1)^{3k}$ ways $\mathcal{A}$ can be selected, so by the union bound
we have
\begin{align*}
&\probnew[ G(n,p) \mbox{ is not $(k,3)$-spoiling }] \leq (n+1)^{3k}\exp\left(-(n-3k)(1-p)^k\right) \\
&\leq \exp\left(3\frac{\ln n-2\ln\ln n+\ln\ln \frac{1}{1-p}-\ln 3}{-\ln(1-p)} \ln n -n\frac{3\ln^2 n}{-n\ln (1-p)}+o(1)\right)
 =o(1),
\end{align*}
assuming $p\leq 1-\frac{1}{n}$.
Otherwise $k=0$ and the statement of the lemma becomes trivial, since every graph is $0$-spoiling.
\proofend

The next lemma provides the lower bound in Theorem~\ref{random-graph}
when $p\leq \frac{1}{2}$.
\begin{lem}
For every $\epsilon>0$ there exists a constant $K=K(\epsilon)$ such
that for all $p$ with
$K/n\leq p\leq 1/2$, the graph $G(n,p)$ is a.a.s.\ $k$-spoiling for
$k=\left\lfloor(1-\epsilon)\frac{\ln(np)}{-\ln(1-\mu )}\right\rfloor$.
\end{lem}

\proofstart
Similarly to the last section, we fix $m=\left\lfloor
\frac{1}{p}\right\rfloor$. We show
that $G(n,p)$ is a.a.s.\ $(k,6m)$-spoiling.
First we observe that for any fixed $S\subset V$
of size at most $6m $ and a fixed vertex $x \in V\setminus S$, the
probability that $x$ spoils $S$ is at least $1-\mu$. Note
that $1-\mu$ is exactly the probability if $|S|=m$ and by the definition of
$\mu$ as a maximum it is at least this much for other sizes strictly below
$6m$.
A simple way to see the bound for $|S|=6m$ is to partition $S$ into six parts
of size $m$ each.
The probability that $x$ has exactly one neighbor in any one of them is $\mu$,
these events are independent, so the probability that this holds for at least
two of them is exactly $1-(1-\mu)^6 -6\mu(1-\mu)^5$. Since we have
$\mu>e^{-1}$ this is larger than $1-\mu$.

Note that $k<  3 \ln (np)$ since $\mu>1/e$.
First we fix a family $\mathcal{A}$ of
$k$ disjoint sets of size at most $6m$ each
and estimate the probability that no vertex $x \in V\setminus
\bigcup {\mathcal A}$ is a spoiler for it.
\begin{align*}
\probnew\left[\nexists x \in V\setminus \bigcup\mathcal{A}: x\mbox{ spoils } \mathcal{A}\right ]
&\leq \left(1-(1-\mu) ^{k}\right)^{n-\left|\bigcup\mathcal{A}\right|}\\
&\leq \exp\left( -\frac{n}{2}(1-\mu)^k\right)\\
&\leq \exp\left(-\frac{(np)^{\epsilon}}{2p}\right),
\end{align*}
where in the second inequality we use the fact that $\left|\bigcup\mathcal A\right|\le6mk<n/2$ for $K$ large enough.

The union bound for the probability that this happens for
any family $\mathcal{A}$ of $k$ sets of size at most $6m$ each is enough now to
finish the proof:
\begin{align*}
\probnew\left[\exists \mathcal A \forall x\in V\setminus \bigcup\mathcal{A}: x\mbox{ does not spoil } \mathcal{A}\right ]
&< \left(\sum_{i\leq 6m}\binom{n}{i}\right)^k \exp\left(-\frac{(np)^{\epsilon}}{2p}\right)\\
&<(np)^{6mk}\exp\left(-\frac{(np)^{\epsilon}}{2p}\right)\\
& < \exp\left(\frac{18\ln^2(np)}{p}-\frac{(np)^{\epsilon}}{2p}\right)\\
&=o(1).
\end{align*}
\proofend

\section{Graphs with large conflict-free chromatic number}
\label{sec:main}

In this section we show the existence of $n$-vertex graphs $G$
with $\chi_{CF}(G) = \Omega(\ln^2 n)$.
This gives the correct order of magnitude of $\chi_{CF} (n)$ and proves
Theorem~\ref{main}.

To show the statement, we construct an $n$-vertex graph $G$ using random methods.
The vertex set is partitioned into classes $L_1, \ldots , L_k$ of size
$\frac{n}{k}$ each, with $k=\lfloor\ln n\rfloor$.
The edges will be selected at random, independently of each other. To define the
probabilities
we let the {\em weight} of a vertex $x\in L_i$ be \[w_x=0.99^{i}.\]
The probability of an edge between vertices $x\in L_i$ and $y\in L_j$ is equal to
\[\probnew[xy\in E(G)] := w_xw_y= 0.99^{i+j}.\]

The {\em weight} of a set $S\subseteq V$ is defined to be the sum of the weights of
its elements,
\[w(S)=\sum_{v\in S}w_v.\]

For a vertex coloring $\chi$ we say that {\em vertex $v$
takes care of itself} if
the color of $v$ is unique in $N[v]$,
i.e. every $u\in N(v)$ has a color different from $\chi(v)$.
We say that {\em a color class $S$ takes care of a vertex $x$}
if $x\in N^{(1)}(S)$.
The crucial probability, denoted by $p(x,S)$, that a vertex $x\in L_i$ is taken care
of by a color class $S$ not containing $x$ is equal to
\begin{align*}
p(x,S)=\probnew[ |N(x)\cap S|=1]
&= \sum_{s\in S}\probnew[N(x)\cap S = \{s\} ] \\
&= \sum_{s\in S}w_sw_x \prod_{y\in S\setminus\{s\}}(1-w_yw_x) \\
&< w_x\sum_{s\in S}w_s \exp\left( - \sum_{y\in S\setminus\{s\}} w_yw_x\right) \\
&= w_x\sum_{s\in S}w_s \exp\left( - w(S)w_x +w_sw_x\right) \\
&\leq w_xw(S) e^{-w_xw(S)+0.99}.
\end{align*}
Note that since the function $ze^{-z}$ has a unique maximum at $z=1$, we
always have $p(x,S) < e^{-0.01}$.
If $\chi$ is a conflict-free coloring,
then every vertex is taken care of either by itself or by a color class not containing this vertex.

We call a set {\em heavy} if its weight is larger than $\sqrt{n}$, otherwise we call
it {\em light}.
Note that since any vertex has weight at least $0.99^{\ln n}> n^{-0.02}$, we obtain for any light color class $S$
\[|S|< w(S) n^{0.02}< n^{0.52}.\]

In the following lemma we list three properties, which hold a.a.s.\
for our random $G$ and, together, imply that no conflict-free coloring exists
with $o\left(\ln^2 n\right)$ colors.
As usual, $\alpha(G)$ denotes the {\em independence number} of $G$,
i.e. the size of a largest independent set.
\begin{lem}\label{threecond}
For $G$ the following three properties hold a.a.s.
\begin{itemize}
\item[$(i)$] $\alpha(G) \leq n^{0.6}$.
\item[$(ii)$] For every heavy set $S\subseteq V$,
we have $\left|N^{(1)} (S)\right|<n^{0.6}$.
\item[$(iii)$] Let $r=\left\lfloor10^{-5}\ln^2 n\right\rfloor$.
For all pairwise disjoint light sets $S_1, \ldots , S_{r}\subseteq V$,
we have $\left|\bigcup_{i=1}^r N^{(1)}(S_i)\right|
<n-n^{0.7}$.
\end{itemize}
\end{lem}
\proofstart
Since the probability for each pair of vertices to be an edge
of $G$ is at least $0.99^{2\ln n}$,
the largest independent set is at most as large as it is in
$G(n, 0.99^{2\ln n})$.
It is well-known (see, for example Theorem 11.25 (ii) in Bollob\'as's book~\cite{bollobas} for more details) that for $2.27/n\leq p\leq 1/2$, a.a.s.\
the largest independent set in $G(n,p)$ has size at most $2\frac{\ln\left(np\right)}{p}$.
Thus, the largest independent set in $G$ a.a.s.\ has size at most $2\frac{\ln\left(0,99^{2\ln n}n\right)}{0,99^{2\ln n}}< n^{0.6}$.

For the second statement fix a subset $S\subseteq V$
with weight at least $n^{0.5}$ and a set $A\subseteq V\setminus S$
with at least $n^{0.6}$ elements.
We estimate the probability that all elements $x\in A$ are in the one-neighborhood of
$S$.
\begin{align*}
\probnew\left[ N^{(1)}(S) \supseteq A\right] & = \prod_{x\in A}p(x,S)\\
& \leq \prod_{x\in A} w_xw(S) e^{-w_xw(S) + 1} \\
&< \left(n^{0.48} \exp \left(-n^{0.48} +1\right)\right)^{n^{0.6}} \\
& = \exp\left( -n^{1.08}(1+o(1))\right).
\end{align*}
(Here we used that $w_xw(S)> 0.99^{\ln n}n^{0.5} > n^{0.48}$
and that $ze^{-z}$ is decreasing in the interval $[1, \infty)$.)
Summing up over all the at most $2^n\cdot 2^n$ choices of $S$ and $A$
we obtain that the probability that $(ii)$ fails tends to $0$.

For the third part
fix subsets $S_1, \ldots , S_r$ with $w(S_i)\leq \sqrt{n}$
and $B$ with
$|B|=n^{0.7}$.
We estimate the probability that all $x\in V\setminus B$ are
in the one-neighborhood of at least one of the $S_i$.

For this we first show that $\sum_{i=1}^r p(x,S_i) > 0.01 \ln n$
for at most half of the vertices $x\in V$.
Indeed, otherwise
\begin{align*}
\frac{n}{2}\cdot 0.01\ln n
& \leq \sum_{x\in V} \sum_{i=1}^r p(x,S_i) \\
& = \sum_{i=1}^r \sum_{x\in V} p(x,S_i) \\
& \leq r \left(100e+100+200\right) \frac{n}{\ln n},
\end{align*}
contradicting the definition of $r$.
For the last estimate we used that for a fixed color class $S_i$,
$\sum_{x\in V} p(x,S_i) \leq \frac{n}{\ln n}\sum_{j=1}^\infty z_j e^{-z_j+1}$,
where $z_j$ is a geometric progression with quotient $0.99$.
The terms of the sum for $z_j\leq 1$ can be estimated by $ez_j$ and hence this
part is at most $\frac{e}{1-0.99}=100e$.
The sum of the terms for
$z_j\geq 2$ can be estimated by
$100\int_1^{\infty} z e^{-z+1}dz =200$.
And finally the sum of the terms for
$1<z_j< 2$ can be estimated by
$100$, since there are at most $100$ such $z_j$'s, and for each of them the value of the function is at most $1$.

Let $V'\subseteq V$ be the set of those vertices $x\in V$
for which $\sum_{i=1}^r p(x,S_i) \leq 0.01 \ln n$. Then by the above
$|V'| \geq n/2$.
\begin{align*}
\probnew[\forall x\in V\setminus B\ \exists i \mbox{ with } |N(x)\cap S_i|=1]
& = \prod_{x\in V\setminus B} \left(1-\prod_{i=1}^r(1-p(x,S_i))\right) \\
& \leq \exp\left( -\sum_{x\in V'\setminus B} \prod_{i=1}^r(1-p(x,S_i))\right) \\
& \leq \exp\left( -\sum_{x\in V'\setminus B} e^{-5\sum_{i=1}^r p(x,S_i)}\right) \\
& \leq \exp\left( - \left(\frac{n}{2}-n^{0.7}\right) e^{-0.05 \ln n}\right)\\
& \leq \exp\left( -n^{0.95}(1/2-o(1))\right)
\end{align*}
Here we used that in the range of our interest, i.e. for
$0<z=p(x,S_i) < e^{-0.01}$, we have $1-z > e^{-5z}$.

The sets $S_1, \ldots S_r$ and $B$ with the given properties can be
chosen at most
\[{n\choose n^{0.7}}\left((n+1)^{\sqrt{n}}\right)^r =
e^{O\left(n^{0.7}\ln n\right)}\]
ways,
where we first choose a set $B$ of size $n^{0.7}$ from $V$,
and then choose one by one the vertices forming the sets $S_1,\ldots, S_r$.
Hence with probability tending to $1$ the third condition holds.
\proofend

Finally, we show how the above properties imply the existence of graphs without a
conflict-free coloring with $10^{-5}\ln^2 n$ colors.

\textit{Proof of Theorem~\ref{main}.}
Let us take a graph $G$ having properties $(i)-(iii)$ of Lemma~\ref{threecond}
with a sufficiently large vertex set.
Take an arbitrary $r$-coloring $c$ of $G$, where $r=\left\lfloor10^{-5}\ln^2
n\right\rfloor$ as in the lemma.
We prove that $c$
is not a conflict-free coloring.
We define the following sets:\\[2mm]
$\bullet$ The set of all vertices that take care of
themselves, that is \\[1mm]
\phantom{$\bullet$} $T=\{x\in V:~\forall y\in N(x), c(x)\neq c(y)\}$\\[2mm]
$\bullet$ The set of all vertices that are taken care of by a heavy color class, that is \\[1mm]
\phantom{$\bullet$} $H=\left\{x\in V:~\exists z\in N(x): \left|c^{-1}(c(z))\right|> \sqrt{n}\wedge \forall y\in N[x]\setminus\{z\}, c(z)\neq c(y)\right\}$\\[2mm]
$\bullet$ The set of all vertices that are taken care of by a light color class, that is\\[1mm]
\phantom{$\bullet$} $L=\left\{x\in V:~\exists z\in N(x): \left|c^{-1}(c(z))\right|\leq \sqrt{n}\wedge \forall y\in N[x]\setminus\{y\}, c(z)\neq c(y)\right\}.$\\[2mm]
If $c$ were a conflict-free coloring, then $V=T\cup H\cup L$.

{\em Vertices taking care of themselves.}
A set of vertices that take care of themselves and have the
same color must form an independent set in $G$.
Hence by $(i)$ any color class can contain at most
$n^{0.6}$ vertices that take care of themselves. So
$|T|\leq rn^{0.6}$.

{\em Vertices taken care of by heavy color classes.}
Fix a heavy color class $S$. By $(ii)$ at most $n^{0.6}$
vertices are taken care of by $S$.
Hence $|H|\leq rn^{0.6}$.

{\em Vertices taken care of by light color classes.}
Let $S_1, \ldots, S_{r^*}$ be the light color classes of $\chi$.
 By $(iii)$ at most $n-n^{0.7}$
vertices are taken care of by the $S_i$'s.
Hence $|L|\leq n-n^{0.7}$.

Thus, $|T\cup H\cup L|<n=|V|$, and $c$ is not conflict-free.
This concludes the proof that $c$ is not a conflict-free coloring.
\proofend

\section{Remarks and open problems}

\paragraph{Radio networks.}
Recently, Noga Alon pointed out to us that the lower bound from Theorem~\ref{main} is similar to the one obtained by
Alon, Bar-Noy, Linial, and Peleg~\cite{overlooked}. Here we discuss briefly the relation of the two results.
The notation in \cite{overlooked} is different from the one we use here,
as they look at a much more applied problem. Creating a small dictionary between their and our notations,
they speak about {\em processors} when we have vertices, a {\em radio network} is what we call a graph,
{\em transmitting at step} $i$ corresponds to having color $i$, and the {\em transmission} itself is a color class.
The problem they analyze is the following. At the beginning, one processor (the {\em sender})
has a message $M$, and the process stops when $M$ is delivered to every processor of the network.
The communication in the network works as follows: in step $i$, every processor from transmission $T_i$
that already received $M$ sends it to all adjacent processors. A
processor {\em receives} a message in a given step if
precisely one of its neighbors transmits in this step. If none of its neighbors
transmits, it hears nothing. If more than one neighbor transmits,
a collision occurs and the processor hears only noise.
A sequence of transmissions is a {\em broadcast schedule}
for the sender $s$ in a network
if after applying the transmissions, every processor in the network has a copy of the message.
Alon, Bar-Noy, Linial, and Peleg~\cite{overlooked} showed that the shortest
length of a broadcast schedule is $\Omega(\ln^2 n)$ for some radio networks
with $n$ processors. The matching upper bound of $O(\ln^2n)$ for any radio
networks with $n$ processors was established earlier Bar-Yehuda, Goldreich,
and Itai~\cite{BGI}.

On the one hand, Alon et al.~\cite{overlooked}
do not restrict a processor to be part of only one transmission,
while in our setting, a vertex has exactly one color.
On the other hand, in our setting we do not have any scheduling structure, and a vertex does not have to wait until it receives
the message to ``transmit''.
Hence, none of the lower bounds implies the other immediately.
However, observe that in the proof of Theorem~\ref{main}
we do {\em not} use the fact that color classes are disjoint.
Consequently our construction gives a common generalization of Theorem~2.1
of \cite{overlooked} and our
Theorem~\ref{main}, as follows.

Let us denote by $\chi'_{CF}(n)$ the smallest integer
such that every graph $G$ on $n$ vertices satisfies the following.
There exists a family $\mathcal{F}\subseteq 2^{V(G)}$ of subsets of $V(G)$ of size $|\mathcal{F}|\leq \chi'_{CF}(n)$
such that for every vertex $x\in V(G)$,
there exists a set $F\in \mathcal{F}$ with $|N_G(x)\cap F|=1$.
Both our paper and \cite{overlooked} deal with problems with further
requirements on the family $\mathcal F$. We insist that they form a partition
of the vertex set, whereas Alon et al. require an
ordering of the sets in $\mathcal F$ with certain properties to exist.
Proving a lower bound for $\chi'_{CF}$
therefore implies the corresponding lower bounds in both papers.
The construction in \cite{overlooked} has chromatic number $2$,
so it does not provide a meaningful lower
bound for $ \chi'_{CF} \leq \chi_{CF} \leq \chi(G)$.
The proof of
Theorem~\ref{main} does work in this more general scenario
and shows a lower bound of order $\ln^2 n$. The
corresponding upper bound follows either from
\cite{PT09} or \cite{BGI}.
\begin{teo}
$\chi'_{CF}(n)=\Theta(\ln^2 n)$.
\end{teo}

\paragraph{Open problems.}
At the two extreme values of $p$ the trivial upper bounds given by
the chromatic number and the domination number plus one
are tight.
For the very sparse range of $p=o(1/n)$ the random graph $G(n,p)$ is a.a.s.\ a tree, hence
both $\chi(G(n,p))$ and $\chi_{CF}(G(n,p))$ are a.a.s.\ $2$. On the other end for $p\geq \frac{1}{2}$
we showed that $|\chi_{CF}-D| \leq 3$.
The particular questions which remain to be answered:

- In what range is $\chi_{CF}(G(n,p))=D(G(n,p))+1$ a.a.s.?

- In what range $\chi(G(n,p))=\chi_{CF}(G(n,p))$ a.a.s.? In particular we would be interested in where the threshold of
$3$-conflict-free colorability is and how much it is different, if at all,
from the threshold of $3$-colorability.

- Does $\chi_{CF} (G(n,p))$ behave in a unimodal way?
For example one might consider the median function and ask whether it is unimodal.

It is an interesting general question to characterize those graphs
where equality holds for $\chi_{CF}(G)=\chi(G)$
 or $\chi_{CF}(G)=D(G)+1$.

By the concentration results of~\cite{WG01,ichanitatibor}
we have the concentration of $\chi_{CF}(G(n,p))$ on two values a.a.s.\ whenever
$\chi_{CF}(G(n,p))=D(G(n,p))+1$. For what range of $p$ does the two-values-concentration hold a.a.s.?
We have a concentration on three values a.a.s.\ whenever $\ln3/\ln(1-p) \approx 0$.
In the worst case, when $p=1/2$, we have concentration on $5$ values a.a.s.
For $p \geq \frac{\sqrt{5}-1}{2}$, we have concentration on four values a.a.s.
(For this we need to consider the $(k,2)$-spoiling property and adapt
the proof of Lemma~\ref{concentration}.)
It would be interesting to obtain a concentration on $2$ values for a wider
range of $p$.

\bigskip

{\bf Acknowledgement.}
We would like to thank Yury Person for pointing out an
error in an earlier version of the proof of Theorem~\ref{main} and Fabian Rabe for useful comments.
Furthermore, we are grateful to Noga Alon for drawing our attention on the
result on radio networks~\cite{overlooked}, as well as for fruitful discussions related to it.


\begin{thebibliography}{}

\bibitem{overlooked}
N. Alon, A. Bar-Noy, N. Linial, and D. Peleg,
A lower bound for radio broadcast,
 {\em Journal of Computer and System Sciences} \textbf{43} (1991), 290--298.

\bibitem{BGI}
R. Bar-Yehuda, O. Goldreich, and A. Itai, On the time-complexity of broadcast in radio
networks: An exponential gap between determinism and randomization, in {\em Proceedings of
the 4th ACM Symposium on Principles of Distributed Computing} (1986), 98--107.

\bibitem{bollobas} B. Bollob\'as, {\bf Random Graphs}, Academic Press, New
York, 2nd ed. (2001).

\bibitem{even} G. Even, Z. Lotker, D. Ron, and S. Smorodinsky, Conflict-free colorings of simple geometric regions
with applications to frequency assignment in cellular networks,
{\em SIAM Journal on Computing}, \textbf{33(1)} (2003), 94--136.

\bibitem{ichanitatibor} R. Glebov, A. Liebenau, and T. Szab\'o,
Concentration of the domination number of the random graph, {\em
  submitted}. preprint available on \url{http://arxiv.org/pdf/1209.3115v3.pdf}.

\bibitem{Luczak} T. \L uczak, Size and connectivity of the $k$-core of a random graph,
{\em Discrete Mathematics} \textbf{91} (1991), 61--68.

\bibitem{PT09} J. Pach and G. Tardos, Conflict-free colorings of graphs and hypergraphs,
{\em Combinatorics, Probability and Computing} \textbf{18} (2009), 819--834.

\bibitem{PSW} B. Pittel, J. Spencer, and N.C. Wormald, Sudden emergence of a giant $k$-core in a random graph,
{\em Journal of Combinatorial Theory, Series B} \textbf{67} (1996), 111--151.

\bibitem{shakharsdiss} S. Smorodinsky, Combinatorial Problems in Computational Geometry, Ph.D Dissertation,
School of Computer Science, Tel-Aviv University (2003).

\bibitem{shakhar} S. Smorodinsky, Conflict-Free Coloring and its Applications,
to appear in
{\em Geometry -- Intuitive, Discrete, and Convex}, (I. Barany, K.J. Boroczky, G. Fejes Toth, and J. Pach, eds.) Bolyai Society Mathematical Studies, Springer.
Preprint available on \url{http://arxiv.org/pdf/1005.3616.pdf}.

\bibitem{WG01} B. Wieland and A.P. Godbole, On the Domination Number of a Random Graph,
{\em Electronic Journal of Combinatorics} \textbf{8(1)} (2001), \#R37.

\end{thebibliography}
\end{document}